\documentclass{article}
\def\Sec{ \mbox{Sec} }

\def\A{{ {\bf A} }}

\def\Sym{ \mbox{Sym} }

\def\1ox{{ \Omega^1_{\scriptstyle{X}} }}
\def\2ox{{ \Omega^2_{\scriptstyle{X}} }}

\def\ok1{{ \Omega^1_K }}
\def\ok2{{ \Omega^2_K }}

\def\P{{ {\bf P} }}

\def\ra{{ \rightarrow }}

\def\Q{{ {\bf Q} }}
\def\D{{ \Delta }}

\def\8{{ {\infty } }}

\def\^{{ ^{\wedge} }}

\def\sup{{ \mbox{sup} }}

\newtheorem{thm}{Theorem}

\newtheorem{lem}{Lemma}

\usepackage{amsmath}
\usepackage{amsfonts}
\usepackage{amscd}
\usepackage{amssymb}

\def\D{{ \Delta }}

\def\D{{ \Delta }}

\title{Relating decision and search algorithms for
 rational points on curves of higher genus}
\author{Minhyong Kim}
\begin{document}
\maketitle
\section{Introduction}
Let 
$$F(x,y)=0 \ \ \ \ \ \ \ \ \ \ (*)$$
be a polynomial equation over the rational numbers of degree $d$
and assume that the complex solutions carve out
an irreducible algebraic curve of geometric genus at least two. This condition
 will be satisfied by
the `generic' equation of degree $\geq 4$.

A theorem of Faltings \cite{Fa} says that the equation has at most
finitely
many rational solutions. The nature of Faltings'
proof, however,
as well as the several other proofs that followed \cite{Vo}, \cite{Bo},
\cite{MW}, 
does not  provide, even in principle,
 an effective
algorithm for actually finding all the solutions.
This is because the proofs provide only a bound on 
 size differences (or ratios), at least in principle,
 rather than
on the absolute size of solutions. (Here, the size of a rational
number can
be taken as the larger of the absolute values of the
numerator and the
denominator.) The {\em effective Mordell
conjecture } proposes one program whereby this
deficiency can be remedied.
Given a rational point in the plane
$P=(p/r,q/r)$ written in reduced form so that
$(p,q,r)=1$, define its {\em height}
to be 
$$h(P):=\sup \{ |p|, |q|, |r|\}$$
Also, given the equation (*), define
$$M(F):=\mbox{sup}\{h(P) \ : \ 
f(P)=0 \},$$
the {\em maximum height} of a rational solution.
The conjecture, in weakest form,
 says
 $M(F)$ should be computable
as a function of $F$. More precise forms of the conjecture
usually
postulate arithmetic-geometric quantities 
that should appear in an effective bound for the maximum size \cite{La2}.

At present, the conjecture appears to be very much 
out of reach and much effort is  expended in relating
one effective conjecture to another or in proving function field
analogues. Because of the lack of palpable structure
on curves of higer genus, as opposed to elliptic
curves, for example, it seems hard to come by
techniques for producing solutions at all, let
alone a complete set. We should note, however,
that much progress has been made in the arena of 
effective algorithms for {\em integral} solutions of
special families of equations, mostly building on
Baker's method \cite{Sm}.

As mentioned, the main application of the effective Mordell
conjecture would be  a search algorithm which would
allows us to deterministically find all solutions
to higher genus equations.
On the other hand, it is conceivable that a search algorithm
exists even without a priori bounds of the sort usually
considered in Diophantine geometry. Put differently,
an algorithm need not necessarily be that given by an
elementary `formula' in terms of geometric invariants,
which is what occurs in the function field case, and
therefore, informs research over number fields.

Even as far as a priori bounds are concerned, consider
the following classical example of genus zero:
$$ax^2+by^2=c$$
with $a,b,c$ positive. A theorem of Holzer \cite{Mo} says that
if a solution exists at all, then there
is a solution $(p/r, q/r)$ in reduced form
that satisfies the bound
$$\sup \{|p|,|q|,|r|\} \leq \sup \{\sqrt{|ab|}, \sqrt{|ac|},\sqrt{|bc| }\}$$
Phrased differently, this is a bound for the {\em minimum}
size of solutions. In this case, since the
solution set can be infinite, there is no question of
a bound on the maximum.
But the utility of the bound shown is that it provides
\medskip

(1) a decision algorithm for the existence of solutions; and
\medskip

(2) the starting point for a {\em purely geometric} algorithm
that allow us to generate all solutions in a precise sense,
even when there are infinitely many.
\medskip

Conjectures of Lang \cite{La3} bounding the  size of  generators
for the Mordell-Weil group of an elliptic curve
are of similar nature.

The purpose of this paper is to express the suspicion 
that this inequality might be the correct paradigm to
follow. That is, perhaps one should formulate conjectures
on the minimum size of the solutions, even in the
higher genus case, and in as elementary a form
as possible. The problem of finding all solutions
could then indeed be
approached in a genuinely algorithmic fashion, rather than by seeking
a simple (or even complicated) formula bounding the maximal size.
Rather than investigate the plausiblity of such conjectures
here (cf. \cite{Sz}), we would like to provide motivation, by 
proving a theorem that  illustrates the utiltiy 
of this point of view in very elementary terms.
Define
$$m(F):=\mbox{inf}\{h(P)\  : \ 
f(P)=0 \},$$
the {\em minimum height} of solutions to $(*)$.

\begin{thm} If $m(F)$ is computable, then
$M(F)$ is computable.
\end{thm}
This theorem is plainly equivalent to
\begin{thm}
A decision algorithm for the existence of a
rational solution to (*) actually provides a determininstic
 search algorithm
for the full solution set.
\end{thm}
{\bf Remarks:}

(1) If
one has a strong enough decision algorithm, one easily
has a search algorithm. For example, if one had a decision
algorithm for  {\em any} algebraic curve presented
in whatever form, 
then simply removing solutions
after we find them through an exhaustive search and
applying the decision algorithm again would tell us
whether we need to continue. To be more concrete,
consider the case where we have a decision algorithm for systems of
equations and {\em inequations}. Then one could construct
a search
 algorithm as follows. Apply the decision algorithm to (*). If the algorithm
says YES, search exhaustively for a solution by ordering the
rational points in $\A^2$ by size. After finding a solution
$(p,q)$ look for all solutions with $x=p$ (easy). Now apply the
algorithm to
$$F(x,y)=0, \ \ \ \ x\neq p$$
to see if we should stop.

Consider also the case where we have an algorithm for
systems of equations in three or more variables.
Then the system above could be replaced by
$$F(x,y)=0, \ \ \ \ \ (x-p)z=1$$
The discerning reader will notice that, in fact,
  a decision algorithm
for systems of two equations in $\A^3$ 
easily  provides a search algorithm. That is,
 a trivial projection technique circumvents the necessity of
 going to progressively higher-dimensional
spaces. 

The point of our
theorem is that it postulates  a decision algorithm
only for  the simplest possible case of a single equation
in two variables. The proof still uses a projection
technique, but of a slightly more subtle nature (Hilbert irreducibility).

In this regard, notice that
Matiyasevich's theorem \cite{Ma}, for example,
 makes it generally desirable to decrease
the number of variables when postulating the existence
of a decision algorithm.

(2) Of course, the kind of decision algorithm we have in
mind is the one mentioned: a bound for the
minimum size of solutions. But from the viewpoint
of recursion theory or computer science, one could imagine
other strategies. The motivaton for this theorem
 arose exactly from the notion that the problem
of finding a decision
algorithm for a single equation in two-variables
has a decidedly elementary flavour that suggests
approaches from  disciplines other than Diophantine
geometry.

\medskip

\section{Proof of theorem.}

The idea is to send rational points off to infinity
as we find them. Assume we are given a decision
algorithm that answers
YES or NO, given the equation.
Apply it to our equation. If the algorithm
says NO, stop. If it says YES, search exhaustively
until we find one, after ordering all the rational points
in the plane by size, for example.
Eventually, we will find a solutions $(p,q)$.
First investigate
for any other solutions with $x=p$ by solving the
single variable equation
$$F(p,y)=0$$

Now, construct a system of equations:
$$F(x,y)=0; \ \ \ (x-p)z=1$$
Clearly, the solutions $(x,y,z)$ of this system are in 1-1
correspondence with the solutions $(x,y)$ of the
original equation for which $x\neq p$. 
We will project this curve back into the plane without
creating any new rational points.
That is, look at the corresponding equations in
$\P^3$:
$$w^dF(x/w,y/w)=0;\ \ \ (x-pw)z=w^2$$
and call the curve it defines $C$.

\begin{lem}
There exists a rational point $m$ in the plane $H:\{w=0\}$
such that projecting from $m$  
maps $C$ birationally to a curve $D$ in $\P^2$ with the property
that any rational point
on  $D$ comes from a rational point on $C$.
\end{lem}

{\em Proof of Lemma.}
If we project from a  rational point $m \in H$, a point
in $D$ will be rational iff the line through
$m$ that corresponds to it is rational. If this line meets
the curve at just one point, then the rational
point in $D$ will have
come from a rational point on our original curve $C$. So
the only way new rational points can be created is
if there is a rational secant line to $C$
passing through $m$.
Thus, we need to find a point $m\in H$ with the property
that all the lines through $m$ that are secant to
$C$ are irrational.
Here's how we achieve this:

Consider the rational secant map
$$[C\times C -\D] \times \P^1 \ra \P^3$$
which sends $(a,b, (t_0:t_1)) $ to the point
$t_0a+t_1b$. That is, it is sent to the corresponding
point on the parametrized
secant line through $a$ and $b$. Suppose
the image is of dimension $\leq 2$. The image
cannot be the plane $H$, because if it were, then
all secant lines would lie in $H$ and hence, $C$
itself would have to lie in $H$ (since $w$ would have infinitely
many zeros on $C$), which it doesn't. Hence, the image meets
$H$ properly, so we will have plenty of rational
points on $H$ which lie on no secant line at all.
Therefore, we need only consider the case where the image has
dimension 3. In this case, the map is generically finite
to one. Also, the locus of the
points in the domain where it's not finite-to-one
has dimension $\leq 2$, so its image has dimension
$\leq 1$. In particular, this image meets $H$ at most in a curve.
Let $Z$ be the portion of $[C\times C-\D]\times \P^1$ lying above $H$. 
 $Z$ consists of triples
$(a,b, (t_0:t_1))$ such that $t_0a+t_1b \in H$.
But a generic  line in $\P^3$ meets a plane in exactly one point
and only finitely many secants lie on $H$,
so given most $(a,b)$, $(t_0:t_1)$ is uniquely determined.
Therefore, the projection to $(a,b)$ defines a birational
map $Z \ra C\times C $. Now when viewed on
$C\times C $, our rational map clearly factors through
the symmetric product
$\Sym^2(C) $ which, in turn, maps
birationally to $\Sec(C)$, the secant variety to $C$.

Thus, we end up with a rational map from
the secant variety: $ \Sec(C)
 \ra H$ which is generically
finite-to-one. Clearly, $\Sym^2(C)$ is irreducible. Also, the
symmetric product of a  curve of genus $\geq 1$ is not birational
to $\P^2$ (it has non-zero global differential
forms),
so the map has degree $>1$.  Recall Hilbert's
irreduciblity theorem \cite{La}, which says that
given a map defined over $\Q$
$$f:X \ra U$$
where $U$ is a non-empty open subset of
$\P^n$ and $X$ is an irreducible variety of dimension $n$, there exists a rational point of $U$
which does not `lift' to a rational point of $X$,
i.e., does not lie in $f(X(\Q))$.
Therefore,
there exists a rational point in
$H$ which does not lie in the image of any rational
point in $\Sec(C) $ 
and with finite inverse image. Thus,
projecting from the point will give a generically
1-1 map on $C$, and the points where its not
1-1 will not map to rational points, proving the lemma.
\bigskip

Now, search exhaustively for the point on $H$
whose existence is guaranteed by the lemma.
One does this, for example on the affine coordinate 
chart where $z\neq 0$, by examining for each point $(0,a,b) \in H$ 
 the set :
$$
\{(w,x,y,w',x',y',t) \in \A^3 \times \A^3 \times \A^1| $$
$$F(x,y)=0, \ \ x-wp=w, \ \ F(x',y')=0,\ \ x'-w'p=w' \ \ $$
$$w+tw'=0,\ \  x+tx'=a, \ \ y+ty'=b\}$$
The previous discussion says that  this equation for
seven unknowns in seven variables has finitely many
solution for generic $(a,b)$. In fact, there is a Groebner
basis algorithm for computing the dimension of the
set for any fixed
$(a,b)$ \cite{CLO}. If the application of this algorithm gives you
dimension one move on to the next point. Whenever,
the dimension turns out to be zero, find all the
solutions and check for rationality. For this,
recall that there is a Groebner basis algorithm
which finds all the solutions to a zero-dimensional
equation by elimination theory \cite{CLO}.
Also, note that rationality of the secant
line  means that two points of intersection with $C$,
$(w,x,y)$ and $(w',x',y')$ are either both rational
or quadratically conjugate. This can be readily checked.
 If there is a rational secant
for $(0,a,b)$ move on to the next point. Eventually and algorithmically
one finds the point $m$ whose existence is guaranteed
by the lemma.
Now, project from $m$ and compute an equation for
the image using elimination theory. Now dehomogenize
by putting the image of $H$ (which is a line) at
infinity.
We have thereby arrived at another affine equation
$F_1(x,y)=0$
whose solution set has cardinality strictly
less than that of $(*)$. Since this is
birational to the original curve,
it still has genus $\geq 2$.
Now apply the
decision algorithm again and iterate the process
above if it says YES.
\bigskip

{\bf Acknowledgements:} I am grateful to Avi Wigderson
for inviting me to participate in the IAS-Park City
summer institute on computational complexity theory.
It was the stimulating atmosphere of this program
that first suggested to me  that the theorem was
plausible. Thanks also to Dinesh Thakur, Alexander Perlis,
and Anand Pillay
for
several useful comments on an earlier verion
of this paper.

This research was supported in part by a grant
from the National Science Foundation.

{\footnotesize DEPARTMENT OF MATHEMATICS, UNIVERSITY OF ARIZONA, TUCSON, AZ 85721, U.S.A.,
and KOREA INSTITUTE FOR ADVANCED STUDY, 207-43 CHEONGRYANGRI-DONG
DONGDAEMUN-GU, SEOUL, KOREA 130-012. E-MAIL: kim@math.arizona.edu}

\end{document}